\documentclass[a4paper,10pt]{amsart}
\usepackage[english]{babel}
\usepackage{amsmath,tikz-cd}
\usepackage{amssymb}
\usepackage{mathrsfs}
\usepackage{enumerate}
\usepackage{mathtools,yfonts}
\usepackage{tikz,tkz-euclide}
\usepackage{pgfplots}
\usepackage{float}
\usepackage{xcolor} 
\definecolor{DarkRed}{RGB}{173,0,0}
\definecolor{LightRed}{RGB}{201,0,0}
\usepackage[
    colorlinks=true,
    linkcolor=DarkRed,
    urlcolor=LightRed,
    citecolor=LightRed
    ]{hyperref}
\usetikzlibrary{arrows}
\usepackage{makecell}
\usepackage{mathdots}
\usepackage{tcolorbox}
\usepackage[]{algorithm2e}
\usepackage{amsmath,amssymb,amsfonts,mathtools}
\usepackage{array,booktabs,longtable}
\usepackage{geometry}
\usepackage{microtype}
\usepackage{xcolor}
\usepackage{textcomp}
\usepackage{hyperref}
\usepackage{listings}
\newcommand{\Pbb}{\mathbb{P}}

\DeclareMathOperator{\Sym}{Sym}

\newtheorem{thm}{Theorem}[section]

\newtheorem{Proposition}[thm]{Proposition}
\newtheorem{Lemma}[thm]{Lemma}
\newtheorem{Corollary}[thm]{Corollary}

\newtheorem*{thm*}{Theorem}

\theoremstyle{definition}

\newtheorem{Definition}[thm]{Definition}
\newtheorem{Remark}[thm]{Remark}
\newtheorem{Construction}[thm]{Construction}
\newtheorem{Example}[thm]{Example}

\newtheorem*{acks}{Acknowledgments}

\newtheorem*{thmA}{Theorem A}
\newtheorem*{thmB}{Theorem B}

\definecolor{wwwwww}{rgb}{0.4,0.4,0.4}

\DeclareMathOperator{\Diag}{Diag}

\DeclareMathOperator{\expdim}{expdim}

\DeclareMathOperator{\Sec}{\mathbb{S}ec}

\DeclareMathOperator{\rank}{rank}

\newcommand\Span[1]{\langle{#1}\rangle}

\hypersetup{pdfpagemode=UseNone}
\hypersetup{pdfstartview=FitH}
		
		\pgfplotsset{compat=1.12}		
\begin{document}

\title{The Alexander-Hirschowitz theorem for neurovarieties}

\author[Alex Massarenti]{Alex Massarenti}
\address{\sc Alex Massarenti\\ Dipartimento di Matematica e Informatica, Universit\`a di Ferrara, Via Machiavelli 30, 44121 Ferrara, Italy}
\email{msslxa@unife.it}

\author[Massimiliano Mella]{Massimiliano Mella}
\address{\sc Massimiliano Mella\\ Dipartimento di Matematica e Informatica, Universit\`a di Ferrara, Via Machiavelli 30, 44121 Ferrara, Italy}
\email{mll@unife.it}

\date{\today}

\begin{abstract}
We study the dimension and identifiability of neurovarieties associated to polynomial neural networks. We give an independent geometric proof that the linear bounds $d_i\geq 2n_i-1$ on the activation degrees imply non defectiveness for any number of outputs, a dimension statement previously obtained from finite identifiability. The proof is based on a direct analysis of the differential of the parameterization. We also investigate secant and Grassmann-secant obstructions outside this range and prove global identifiability for multi-output architectures under the same degree bounds.
\end{abstract}

\maketitle
\tableofcontents
\setcounter{tocdepth}{1}

\section*{Introduction}
The field of Neuroalgebraic Geometry is an emerging area of research
focused on studying function spaces defined by machine learning
models based on algebraic architectures like
Polynomial Neural Networks (PNN).
PNNs represent a distinct class of neural
network architectures where the traditional, typically non-polynomial,
activation functions (like ReLU, sigmoid, or tanh) are replaced with
polynomials.
This substitution gives PNNs a unique theoretical and practical profile, making them a significant area of study in machine learning.

These networks  achieve a competitive experimental performance across
a wide range of tasks. Critically, the polynomial activation
functions naturally facilitate the capture of high-order interactions
between the input features. Unlike simpler activation functions, which
might require a much greater depth or width to implicitly model
complex relationships, PNNs explicitly incorporate multiplicative
terms that represent these intricate,
high-order dependencies. 

This powerful capability has led to their successful employment across
numerous domains ranging from computer vision (tasks like image recognition and object
  detection, \cite{CMB+20,CMB+22,YHN+21}) and  image representation (learning efficient and meaningful
    feature encodings for visual data,\cite{YBJ+22}), to physics (solving complex differential equations or
      modeling physical systems, \cite{BK21}) and finance (applications such as time series prediction or risk modeling, \cite{NM18}).

On the theoretical side, the use of polynomials provides a remarkable
benefit: it allows for a fine-grained theoretical analysis of the
network's properties. The set of functions that a PNN can compute, the
function space associated with the architecture, possesses a specific
geometric structure. These function spaces are frequently referred to as
neuromanifolds and their Zariski closures are  called neurovarieties.

Since polynomials are the building blocks, tools from algebraic
geometry can be rigorously applied to analyze these
neurovarieties. Investigating the properties of such spaces,
particularly their dimension, yields crucial insights into the
network's behavior: 
The dimension and structure of the neurovariety shed light on the
impact of a PNN's architectural choices (such as the layer widths and
the activation degrees of the polynomials) on the overall
expressivity (the set of functions it can approximate) of various
PNN subtypes, including standard feedforward, convolutional,
and self-attention architectures.

All these opportunities made Neuroalgebraic Geometry  an emerging area
of research that attracted many contributions in the last years \cite{KTB19, BT20, Xi21, BBCV21, Sh23, LYPL21, SMK24, MSMTK25, KLW24b,HMK24,SMK25}.

This paper investigates the algebraic structure underlying polynomial
neural networks and their neurovarieties
$\mathcal{V}_{\mathbf n,\mathbf d}$. Aimed also at a computer science audience
interested in the theoretical capacity and complexity of network
architectures, this work connects deep learning structures to
fundamental results in algebraic geometry. 

We focus on a class of feed-forward networks defined by a width vector
${\mathbf n} = (n_0, \dots, n_L)$ and a set of activation exponents
${\mathbf d} = (d_1,
\dots, d_{L-1})$. These networks utilize weight matrices $W_i$ and
activation functions $\sigma_i$ which raise inputs element-wise to
the power $d_i$. The output function $F$ of this architecture is a set
of $n_L$
homogeneous polynomials of total degree $d = \prod_{i=1}^{L-1} d_i$. 

The neurovariety $\mathcal{V}_{\mathbf n,\mathbf
  d}$ is  the Zariski closure of
the image of the map that takes the network
weights $(W_1, \dots, W_L)$ to the projective space of output
polynomial coefficients. Essentially, $\mathcal{V}_{\mathbf n,\mathbf
  d}$ captures all possible
output polynomials that can be realized by the given architecture and
their limits.
A key mathematical challenge is determining the dimension of this
variety, which represents the actual complexity or expressivity of the
network in terms of the number of independent parameters required to
generate the output space. Indeed, the dimension reflects the intrinsic
degrees of freedom of the model and is a simple measure of 
expressivity that is more precise than the raw parameter count.

The dimension of neurovarieties has recently been studied in various papers, for instance
\cite{KTB19, MSMTK25, FRW25}. 
The neurovariety has an expected dimension
$\text{expdim}(\mathcal{V}_{\mathbf n,\mathbf d})$, typically calculated based on the total
number of parameters adjusted for obvious representation ambiguities. If the actual dimension
 is strictly less than this expected dimension, the
neurovariety is called defective, also indicating redundant
parameterization. 

From the viewpoint of learning theory and practice, the actual dimension 
\(\dim\mathcal V_{\mathbf n,\mathbf d}\) is the relevant capacity measure of a Polynomial Neural Network (PNN), sharper than raw parameter count. 
Knowing \(\dim\mathcal V_{\mathbf n,\mathbf d}\) has several consequences:
\begin{itemize}
\item[-] \emph{Model selection and redundancy.} If \(\mathcal V_{\mathbf n,\mathbf d}\) is defective, many weights are intrinsically redundant (flat directions of the coefficient map), so a slimmer architecture can achieve the same expressivity. Dimension gives a principled target for pruning and width/depth tuning.
\item[-] \emph{Identifiability and interpretability.} When the actual dimension is strictly below the free parameter count, the generic fiber of the parameterization is positive dimensional, obstructing global identifiability of the network; see Definition~\ref{def:ident} and Corollary~\ref{cor:ident}. 
Our non–defectiveness and identifiability results guarantee that, generically, distinct functions correspond to finitely many parameterizations modulo the usual symmetries.
\item[-] \emph{Sample complexity and generalization.} Effective capacity (hence the amount of data needed to avoid overfitting) scales with \(\dim\mathcal V_{\mathbf n,\mathbf d}\) rather than with the bare number of weights. 
Overestimating capacity by parameter count leads to unnecessarily large models and poorer generalization under noise.
\item[-] \emph{Optimization and numerical stability.} Rank drops in the Jacobian of the parametrization of the neuromanifold produce ill–conditioned directions and slow training. 
Working at architectures whose neurovarieties are non–defective improves conditioning and removes spurious flatness due to algebraic degeneracies.
\item[-] \emph{Compression and deployment.} A correct dimension estimate bounds the minimal representation size for a task, guiding compression without loss of accuracy.
\end{itemize}

If \(\mathcal V_{\mathbf n,\mathbf d}\) fills its ambient coefficient space, the model can interpolate essentially any degree polynomial, which maximizes expressivity but weakens inductive bias and makes the network more prone to fitting noise. 
In many applications one prefers structured hypothesis classes, i.e. proper subvarieties that encode algebraic constraints reflecting the task; this improves robustness and denoising and aligns with standard regularization principles (see also Remark~\ref{rem:nsb}). 
Of course, in tasks that genuinely require arbitrary polynomial expressivity (e.g. symbolic regression without prior structure), filling may be desirable. 
In all cases, the optimal regime is non–defective (actual dimension equals the expected one), with filling vs non–filling decided by the application: non–filling is generally preferable in supervised learning with noisy data, while filling is justified only when maximal expressivity is itself the goal.

Identifiability answers whether the parameters and,
consequently, the hidden representations of a Neural Network can be
uniquely determined from its response. This is considered to be a 
key question in Neural Network  theory, \cite{S92, BBM23}.
Identifiability is critical in many aspects \cite{LBLGSB18, GBEK22,
  KWWC24, UDBC25}:  to ensure interpretability in representation
learning, to provably obtain disentangled representations, in the study of causal
models, to understand how the architecture affects the inference process and to
support manipulation or stitching of retrained models and representations.

The linear degree bound used in this paper already appears in
\cite[Corollary~17]{UDBC25}: in the notation of that paper, the inequalities
$r_i\ge2d_i-1$ imply finite identifiability of homogeneous polynomial neural
networks. By \cite[Proposition~10]{UDBC25}, this gives maximal affine
dimension. After projectivizing the $n_L$ outputs independently, one obtains
the dimension statement of Theorem~\ref{exdNR} below. Corollary~17 of
\cite{UDBC25} gives finite identifiability; under the same bounds our final
result gives global identifiability when $n_L\ge2$.

We give an independent geometric proof of this statement. We view it as an
Alexander--Hirschowitz type theorem for neurovarieties. For $L=2$ and one
output, the neurovariety is a classical secant variety of a Veronese variety;
see Remark~\ref{secV}. Our proof is based directly on the differential of the
network parameterization. We first write the differential in blocks
corresponding to the different layers. We then prove a recursive
linear-independence lemma showing that, when $d_i\ge2n_i-1$, the contributions
produced at each level separate in the required way. This determines the
kernel explicitly: at a general parameter point the only variations in the
kernel are the usual rescalings of the hidden neurons together with the
projective rescalings of the outputs.

Secant varieties and Grassmann-secant constructions remain useful for the
geometric interpretation of the network and for detecting obstructions when
the degree bound is not satisfied. In particular, they explain the defective
examples discussed in Section~\ref{sec:nondef}, but they are not used as a
transversality assumption in our differential proof. Under the same degree
bounds, when $n_L\ge2$ we further prove global identifiability; see
Theorem~\ref{th:main}.
\begin{thmA}
\textit{Set an architecture associated to the width vector
${\mathbf n}=(n_0,\ldots,n_L)$, with $L\geq2$, and activation degrees
${\mathbf d}=(d_1,\ldots,d_{L-1})$, over a subfield of the complex field.
Assume that
$$
n_i\ge2\quad\text{for }i=0,\ldots,L-1,
\qquad
 d_i\ge2n_i-1\quad\text{for }i=1,\ldots,L-1.
$$
Then the neurovariety $\mathcal V_{\mathbf n,\mathbf d}$ is non defective and
$
\dim\mathcal V_{\mathbf n,\mathbf d}
=\sum_{i=1}^L n_i(n_{i-1}-1).
$}
\end{thmA}
Outside the degree range of Theorem A we discuss geometric obstructions to
non defectiveness coming from positive-dimensional families of intermediate
secant spaces. One of the examples was found using the techniques of
\cite{DR2026-MFA-PNN}.

For multi-output networks the same degree bound also gives global
identifiability. The required uniqueness of the compatible nested
configuration is proved in Section~\ref{sec:nondef} by applying elementary
linear independence of Veronese points successively from the last hidden
layer to the first one; see Theorem~\ref{th:main}.
\begin{thmB}
\textit{Under the assumptions of Theorem A, if $n_L\ge2$, then the associated
polynomial neural network is globally identifiable.}
\end{thmB}

A Magma library is available in Section~\ref{MagmaNV}.

\begin{acks}
While completing this paper, we benefited from attending the 
\begin{center}
\href{https://www.umu.se/en/department-of-mathematics-and-mathematical-statistics/research/mathematical-foundations-of-artificial-intelligence/geumetric-deep-learning-workshop-2025/}{\emph{GeUmetric Deep Learning Workshop 2025}}
\end{center}
at Ume\aa\ University (Sweden), August 19--21, 2025.
We warmly thank the organizers and all participants for creating a friendly and stimulating environment that helped us finalize this work.
In particular, we are grateful to Kathl\'en Kohn and Vahid Shahverdi
for important suggestions on an early version of the preprint. We thank
Kevin Dao for pointing out an error in a previous version and for several
useful comments on the revision.

The authors are members of GNSAGA (INdAM).
A.~Massarenti was supported by the PRIN 2022 project 20223B5S8L ``\emph{Birational geometry of moduli spaces and special varieties}''. M.~Mella was partially
supported by PRIN 2022 project 2022ZRRL4C ``\emph{Multilinear Algebraic Geometry}''.
\end{acks}

\section{Neurovarieties}
Fix a subfield $K\subset\mathbb C$, a width vector $\mathbf n=(n_0,\dots,n_L)$ with $L\geq1$, and weight matrices $W_i\in K^{n_i\times n_{i-1}}$ for $i=1,\dots,L$. Set $W=(W_1,\dots,W_L)$. Fix integers $d_i\geq2$ and let $\sigma_i:K^{n_i}\rightarrow K^{n_i}$ be the element-wise activation $\sigma_i(x_0,\dots,x_{n_i-1})=(x_0^{d_i},\dots,x_{n_i-1}^{d_i})$ for $i=1,\dots,L-1$. Set
$$
F=W_L\circ\sigma_{L-1}\circ W_{L-1}\circ\dots\circ W_2\circ\sigma_1\circ W_1:K^{n_0}\longrightarrow K^{n_L}.
$$
The coordinate functions $F_1,\dots,F_{n_L}$ are homogeneous polynomials of degree $d=\prod_{i=1}^{L-1}d_i$, that is $F\in(\Sym_dK^{n_0})^{n_L}$. Write $F_\ell=\sum_js_{F_\ell,j}x^{\alpha(j)}$, where the degree $d$ monomials $x^{\alpha(j)}$ are ordered lexicographically. The coefficients $s_{F_\ell,j}$ are polynomials in the entries of the weight matrices and define the map
$$
\begin{array}{cccc}
\phi: & K^{n_1\times n_0}\times\cdots\times K^{n_L\times n_{L-1}} & \longrightarrow & \bigl(\mathbb P(\Sym_dK^{n_0})\bigr)^{n_L}\\
& (W_1,\dots,W_L) & \longmapsto & (F_1,\dots,F_{n_L}).
\end{array}
$$
The \textit{neurovariety} $\mathcal V_{\mathbf n,\mathbf d}$ is the Zariski closure of $\operatorname{Im}(\phi)$ in $\bigl(\mathbb P(\Sym_dK^{n_0})\bigr)^{n_L}$.

We fix the entries of the last column of every $W_i$ to $1$. This gives an affine parameter space with $\sum_{i=1}^Ln_i(n_{i-1}-1)$ coordinates. After dehomogenizing the target we get an affine map $\phi_A$, and the rank of $J\phi_A$ at a general point is the dimension of $\mathcal V_{\mathbf n,\mathbf d}$.

\begin{Definition}
The \textit{expected dimension} of $\mathcal V_{\mathbf n,\mathbf d}$ is
$\expdim(\mathcal V_{\mathbf n,\mathbf d}) = \min\left\lbrace\sum_{i=1}^L n_i(n_{i-1}-1) , N\right\rbrace$
where $N = n_{L} \binom{n_0-1+d}{n_0-1}-n_L$. We will say that $\mathcal V_{\mathbf n,\mathbf d}$ is \textit{defective} if $\dim(\mathcal V_{\mathbf n,\mathbf d}) < \expdim(\mathcal V_{\mathbf n,\mathbf d})$.
\end{Definition}

\begin{Remark}\label{secV}
Assume $L = 2, n_2 = 1$, then $\mathcal V_{\mathbf n,\mathbf d}= \Sec_{n_1}(\mathcal{V}^{n_0-1}_{d_1})$ is the secant variety of $(n_1-1)$-planes that are $n_1$-secant to the Veronese variety $\mathcal{V}^{n_0-1}_{d_1}$: the image of the degree $d_1$ Veronese embedding of $\mathbb{P}^{n_0-1}$. Indeed, write
$$
W_1 = \left(
\begin{array}{ccc}
\alpha_{0,0} & \dots & \alpha_{0,n_0-1} \\ 
\vdots & \ddots & \vdots \\ 
\alpha_{n_1-1,0} & \dots & \alpha_{n_1-1,n_0-1}
\end{array} 
\right),\quad
W_2 = (\beta_0,\dots,\beta_{n_1-1}).
$$
Then
\begin{footnotesize}
$$
\begin{array}{lll}
F & = & (W_2\circ\sigma_1\circ W_1)(x_0,\dots , x_{n_0-1}) = (W_2\circ\sigma_1)(\alpha_{0,0}x_0+\dots + \alpha_{0,n_0-1}x_{n_0-1},\dots ,\alpha_{n_1-1,0}x_0+\dots + \alpha_{n_1-1,n_0-1}x_{n_0-1})) \\ 
 & = & W_2((\alpha_{0,0}x_0+\dots + \alpha_{0,n_0-1}x_{n_0-1})^{d_1},\dots ,(\alpha_{n_1-1,0}x_0+\dots + \alpha_{n_1-1,n_0-1}x_{n_0-1})^{d_1})\\
 & = & \beta_{0}(\alpha_{0,0}x_0+\dots + \alpha_{0,n_0-1}x_{n_0-1})^{d_1} + \dots + \beta_{n_1-1}(\alpha_{n_1-1,0}x_0+\dots + \alpha_{n_1-1,n_0-1}x_{n_0-1})^{d_1}.
\end{array} 
$$
\end{footnotesize}
and hence $F\in \mathbb{P}(\operatorname{Sym}_d K^{n_0})$ lies in
$\Sec_{n_1}(\mathcal{V}^{n_0-1}_{d_1})$. The Alexander-Hirschowitz
theorem \cite[Theorem 1]{AH95} states that when $K={\mathbb C}$ all such secant varieties attain the expected dimension except for the following defective cases:

\begin{center}\renewcommand{\arraystretch}{1.15}
\begin{tabular}{c|c|c}
$n_0$ & $d_1$ & $n_1$ \\ \hline
$n_0\ge 3$ & $2$ & $2\leq n_1 \leq n_0-1$ \\
$3$ & $4$ & $5$\\
$4$ & $4$ & $9$\\
$5$ & $3$ & $7$\\
$5$ & $4$ & $14$\\
\hline
\end{tabular}
\end{center}

Note that if $K\subset{\mathbb C}$ is a subfield it contains the
rational numbers. The equations of the Veronese variety are
rational therefore on any field $K\subset{\mathbb C}$ the secant
varieties of the Veronese variety that are non defective over
${\mathbb C}$ are non defective over $K$. 
\end{Remark}

It is known that in polynomial neural networks several different parameters lead to the same
polynomial. Indeed, for each hidden layer one can permute the hidden neurons, and rescale the
input and output to each activation function, \cite{UDBC25} \cite{KTB19}. These transformations
lead to a different set of parameters that leave the neurovariety
unchanged. This observation motivates the following definition.

\begin{Definition}[\cite{UDBC25}]\label{def:ident}
The polynomial neural network is said to be globally identifiable if for a
general element in the neurovariety all its representations  can be obtained from the permutations and
element-wise scalings of a single representation.
\end{Definition}
In Section~\ref{sec:2} we will reinterpret this definition from a more
geometric point of view that will be easily applied in our context,
see Corollary~\ref{cor:ident}.

\begin{Construction}
Set
$$
W_i = \left(
\begin{array}{ccc}
\alpha_{0,0}^i & \dots & \alpha_{0,n_{i-1}-1}^i \\ 
\vdots & \ddots & \vdots \\ 
\alpha_{n_i-1,0}^i & \dots & \alpha_{n_i-1,n_{i-1}-1}^i
\end{array} 
\right)
$$
for $i = 1,\dots,L$.
\begin{itemize}
\item[($s_0$)] At step zero we have $\mathbb{P}^{n_0-1}$ with homogeneous coordinates $[x_0:\dots : x_{n_0-1}]$.
\item[($s_1$)] We apply $W_1$ producing $n_1$ linear forms 
$L_j = \sum_{i=0}^{n_0-1}\alpha_{j,i}^1x_i, \, j = 0,\dots,n_1-1$.
\item[($s_2$)] We raise the $L_j$ to the $d_1$ power and apply $W_2$ producing $n_2$ polynomials
$F_{2,j} = \sum_{i=0}^{n_1-1}\alpha_{j,i}^2 L_i^{d_1}, \, j = 0,\dots,n_2-1$. 
Note that we may interpret $F_{2,j}$
\begin{itemize}
\item[($i_1$)] either as a polynomial of degree $d_1$ on $\mathbb{P}^{n_0-1}$;
\item[($i_2$)] or as a polynomial of degree $1$ on $\mathbb{P}^{M_1}$ with $M_1 = \binom{n_0-1+d_1}{n_0-1}-1$.
\end{itemize}
However, these two interpretations coincide in this particular step.
\item[($s_3$)] We raise the $F_{2,j}$ to the $d_2$ power and apply $W_3$ producing $n_3$ polynomials
$F_{3,j} = \sum_{i=0}^{n_2-1}\alpha_{j,i}^3 F_{2,i}^{d_2}, \, j = 0,\dots,n_3-1$. 
We may interpret $F_{3,j}$
\begin{itemize}
\item[($i_1$)] either as a polynomial of degree $d_2$ on $\mathbb{P}^{M_1}$;
\item[($i_2$)] or as a polynomial of degree $1$ on $\mathbb{P}^{M_2}$ with $M_2 = \binom{M_1+d_2}{M_1}-1$.
\end{itemize}
From this step onward the two interpretations actually differ.
\item[$\vdots$]
\item[($s_k$)] At step $k$ we raise the $F_{k-1,j}$ to the $d_{k-1}$ power and apply $W_k$ producing $n_k$ polynomials
$F_{k,j} = \sum_{i=0}^{n_{k-1}-1}\alpha_{j,i}^k F_{k-1,i}^{d_{k-1}}, \, j = 0,\dots,n_k-1$, 
where the $F_{k-1,i}$ are the $n_{k-1}$ polynomials produced at step ($s_{k-1}$). We may interpret $F_{k,j}$
\begin{itemize}
\item[($i_1$)] either as a polynomial of degree $d_{k-1}$ on $\mathbb{P}^{M_{k-2}}$;
\item[($i_2$)] or as a polynomial of degree $1$ on $\mathbb{P}^{M_{k-1}}$ with $M_{k-1} = \binom{M_{k-2}+d_{k-1}}{M_{k-2}}-1$.
\end{itemize}
\item[$\vdots$]
\item[($s_L$)] At the final step $L$ we raise the $F_{L-1,j}$ to the $d_{L-1}$ power and apply $W_L$ producing $n_L$ polynomials
$F_{L,j} = \sum_{i=0}^{n_{L-1}-1}\alpha_{j,i}^L F_{L-1,i}^{d_{L-1}}, \, j = 0,\dots,n_L-1$, 
where the $F_{L-1,i}$ are the $n_{L-1}$ polynomials produced at step ($s_{L-1}$). We may interpret $F_{L,j}$
\begin{itemize}
\item[($i_1$)] either as a polynomial of degree $d_{L-1}$ on $\mathbb{P}^{M_{L-2}}$;
\item[($i_2$)] or as a polynomial of degree $1$ on $\mathbb{P}^{M_{L-1}}$ with $M_{L-1} = \binom{M_{L-2}+d_{L-1}}{M_{L-2}}-1$.
\end{itemize}
\end{itemize}
\end{Construction}

\begin{Example}
Consider the case $(n_0,n_1,n_2,n_3) = (2,2,2,1)$, $(d_1,d_2) = (2,2)$. We have 
\begin{tiny}
$$
F = c_{0,0}(b_{0,0}(a_{0,0}x_0 + a_{0,1}x_1)^2 + b_{0,1}(a_{1,0}x_0 + a_{1,1}x_1)^2)^2 + c_{0,1}(b_{1,0}(a_{0,0}x_0 + a_{0,1}x_1)^2 + b_{1,1}(a_{1,0}x_0 + a_{1,1}x_1)^2)^2.
$$
\end{tiny}
The coefficients of $F$ as a polynomial of degree $4$ in $x_0,x_1$ are:
\begin{tiny}
$$
\begin{array}{lll}
s_0 & = & a_{0,0}^4b_{0,0}^2c_{0,0} + a_{0,0}^4b_{1,0}^2c_{0,1} + 2a_{0,0}^2a_{1,0}^2b_{0,0}b_{0,1}c_{0,0} + 2a_{0,0}^2a_{1,0}^2b_{1,0}b_{1,1}c_{0,1} + a_{1,0}^4b_{0,1}^2c_{0,0} + a_{1,0}^4b_{1,1}^2c_{0,1};
 \\ 
s_1 & = & 4a_{0,0}^3a_{0,1}b_{0,0}^2c_{0,0} + 4a_{0,0}^3a_{0,1}b_{1,0}^2c_{0,1} + 4a_{0,0}^2a_{1,0}a_{1,1}b_{0,0}b_{0,1}c_{0,0}  + 4a_{0,0}^2a_{1,0}a_{1,1}b_{1,0}b_{1,1}c_{0,1} + 4a_{0,0}a_{0,1}a_{1,0}^2b_{0,0}b_{0,1}c_{0,0}\\
 & & + 4a_{0,0}a_{0,1}a_{1,0}^2b_{1,0}b_{1,1}c_{0,1} + 4a_{1,0}^3a_{1,1}b_{0,1}^2c_{0,0} + 4a_{1,0}^3a_{1,1}b_{1,1}^2c_{0,1};\\ 
s_2 & = & 6a_{0,0}^2a_{0,1}^2b_{0,0}^2c_{0,0} + 6a_{0,0}^2a_{0,1}^2b_{1,0}^2c_{0,1} + 2a_{0,0}^2a_{1,1}^2b_{0,0}b_{0,1}c_{0,0} + 2a_{0,0}^2a_{1,1}^2b_{1,0}b_{1,1}c_{0,1} + 8a_{0,0}a_{0,1}a_{1,0}a_{1,1}b_{0,0}b_{0,1}c_{0,0} \\
 & & + 8a_{0,0}a_{0,1}a_{1,0}a_{1,1}b_{1,0}b_{1,1}c_{0,1} + 2a_{0,1}^2a_{1,0}^2b_{0,0}b_{0,1}c_{0,0} + 2a_{0,1}^2a_{1,0}^2b_{1,0}b_{1,1}c_{0,1} + 6a_{1,0}^2a_{1,1}^2b_{0,1}^2c_{0,0} + 6a_{1,0}^2a_{1,1}^2b_{1,1}^2c_{0,1};\\ 
s_3 & = & 4a_{0,0}a_{0,1}^3b_{0,0}^2c_{0,0} + 4a_{0,0}a_{0,1}^3b_{1,0}^2c_{0,1} + 4a_{0,0}a_{0,1}a_{1,1}^2b_{0,0}b_{0,1}c_{0,0} + 4a_{0,0}a_{0,1}a_{1,1}^2b_{1,0}b_{1,1}c_{0,1} + 4a_{0,1}^2a_{1,0}a_{1,1}b_{0,0}b_{0,1}c_{0,0}  \\
 & & + 4a_{0,1}^2a_{1,0}a_{1,1}b_{1,0}b_{1,1}c_{0,1} + 4a_{1,0}a_{1,1}^3b_{0,1}^2c_{0,0} + 4a_{1,0}a_{1,1}^3b_{1,1}^2c_{0,1}; \\
s_4 & = & a_{0,1}^4b_{0,0}^2c_{0,0} + a_{0,1}^4b_{1,0}^2c_{0,1} + 2a_{0,1}^2a_{1,1}^2b_{0,0}b_{0,1}c_{0,0} + 2a_{0,1}^2a_{1,1}^2b_{1,0}b_{1,1}c_{0,1} + a_{1,1}^4b_{0,1}^2c_{0,0} + a_{1,1}^4b_{1,1}^2c_{0,1}.
\end{array} 
$$
\end{tiny}
Let $G(z_0,z_1,z_2)$ be the degree $2$ polynomial obtained by setting $z_0 = x_0^2,z_1 = x_0x_1,z_2 = x_1^2$ in $F$. The coefficients of $G$ are:
\begin{tiny}
$$
\begin{array}{lll}
f_0 & = & a_{0,0}^4b_{0,0}^2c_{0,0} + a_{0,0}^4b_{1,0}^2c_{0,1} + 2a_{0,0}^2a_{1,0}^2b_{0,0}b_{0,1}c_{0,0} + 2a_{0,0}^2a_{1,0}^2b_{1,0}b_{1,1}c_{0,1} + a_{1,0}^4b_{0,1}^2c_{0,0} + a_{1,0}^4b_{1,1}^2c_{0,1}; \\ 
f_1 & = & 4a_{0,0}^3a_{0,1}b_{0,0}^2c_{0,0} + 4a_{0,0}^3a_{0,1}b_{1,0}^2c_{0,1} + 4a_{0,0}^2a_{1,0}a_{1,1}b_{0,0}b_{0,1}c_{0,0} + 4a_{0,0}^2a_{1,0}a_{1,1}b_{1,0}b_{1,1}c_{0,1} + 4a_{0,0}a_{0,1}a_{1,0}^2b_{0,0}b_{0,1}c_{0,0}\\ 
 & & + 4a_{0,0}a_{0,1}a_{1,0}^2b_{1,0}b_{1,1}c_{0,1} + 4a_{1,0}^3a_{1,1}b_{0,1}^2c_{0,0} + 4a_{1,0}^3a_{1,1}b_{1,1}^2c_{0,1};\\ 
f_2 & = & 2a_{0,0}^2a_{0,1}^2b_{0,0}^2c_{0,0} + 2a_{0,0}^2a_{0,1}^2b_{1,0}^2c_{0,1} + 2a_{0,0}^2a_{1,1}^2b_{0,0}b_{0,1}c_{0,0} + 2a_{0,0}^2a_{1,1}^2b_{1,0}b_{1,1}c_{0,1} + 2a_{0,1}^2a_{1,0}^2b_{0,0}b_{0,1}c_{0,0} + 2a_{0,1}^2a_{1,0}^2b_{1,0}b_{1,1}c_{0,1}\\
 & &  + 2a_{1,0}^2a_{1,1}^2b_{0,1}^2c_{0,0} + 2a_{1,0}^2a_{1,1}^2b_{1,1}^2c_{0,1};\\ 
f_3 & = & 4a_{0,0}^2a_{0,1}^2b_{0,0}^2c_{0,0} + 4a_{0,0}^2a_{0,1}^2b_{1,0}^2c_{0,1} + 8a_{0,0}a_{0,1}a_{1,0}a_{1,1}b_{0,0}b_{0,1}c_{0,0} + 8a_{0,0}a_{0,1}a_{1,0}a_{1,1}b_{1,0}b_{1,1}c_{0,1} + 4a_{1,0}^2a_{1,1}^2b_{0,1}^2c_{0,0}\\ 
 & & + 4a_{1,0}^2a_{1,1}^2b_{1,1}^2c_{0,1};\\ 
f_4 & = & 4a_{0,0}a_{0,1}^3b_{0,0}^2c_{0,0} + 4a_{0,0}a_{0,1}^3b_{1,0}^2c_{0,1} + 4a_{0,0}a_{0,1}a_{1,1}^2b_{0,0}b_{0,1}c_{0,0} + 4a_{0,0}a_{0,1}a_{1,1}^2b_{1,0}b_{1,1}c_{0,1} + 4a_{0,1}^2a_{1,0}a_{1,1}b_{0,0}b_{0,1}c_{0,0}\\ 
 & & + 4a_{0,1}^2a_{1,0}a_{1,1}b_{1,0}b_{1,1}c_{0,1} + 4a_{1,0}a_{1,1}^3b_{0,1}^2c_{0,0} + 4a_{1,0}a_{1,1}^3b_{1,1}^2c_{0,1};\\
f_5 & = & a_{0,1}^4b_{0,0}^2c_{0,0} + a_{0,1}^4b_{1,0}^2c_{0,1} + 2a_{0,1}^2a_{1,1}^2b_{0,0}b_{0,1}c_{0,0} + 2a_{0,1}^2a_{1,1}^2b_{1,0}b_{1,1}c_{0,1} + a_{1,1}^4b_{0,1}^2c_{0,0} + a_{1,1}^4b_{1,1}^2c_{0,1}.
\end{array} 
$$
\end{tiny}
Setting $a_{0,1}=a_{1,1}=b_{0,1}=b_{1,1}=c_{0,1} = 1$ we get two maps: $s:\mathbb{A}^5\dashrightarrow\mathbb{P}^4$ given by the $s_i$, and $f:\mathbb{A}^5\dashrightarrow\mathbb{P}^5$ given by the $f_i$, fitting in the following diagram:
$$
\begin{tikzcd}
                                                               &  & \mathbb{P}^5 \arrow[d, "\pi", dashed] \\
\mathbb{A}^5 \arrow[rru, "f", dashed] \arrow[rr, "s"', dashed] &  & \mathbb{P}^4                         
\end{tikzcd}
$$
where $\pi$ is the linear projection from the point of $\mathbb{P}^5$ representing the polynomial $z_0z_2-z_1^2$.
\end{Example}

\section{Neural networks with \texorpdfstring{$n_L=1$}{nL=1}}
\label{sec:2}
Assume $n_L=1$. We describe the corresponding neurovariety in terms of iterated Veronese embeddings and linear projections. This description will be used to interpret the dimension and identifiability problems geometrically.

\begin{Construction}
Set $n_L = 1$ and consider the following diagram:
\begin{footnotesize}
$$
\begin{tikzcd}
                                          &                                                                                                                                                &                                                                                                                         &                                                                                                                          &                                                                 &                                                                                                                                        & \mathcal{V}_{L-1}\subset\mathbb{P}^{M_{L-1}} \arrow[ddddd, "\pi_{L-1}", dashed]                       \\
                                          &                                                                                                                                                &                                                                                                                         &                                                                                                                          &                                                                 & \mathcal{V}_{L-2}\subset\mathbb{P}^{M_{L-2}} \arrow[ru, "\nu_{d_{L-1}}"] \arrow[dddd, "\pi_{L-2}", dashed]                                 &                                                                                           \\
                                          &                                                                                                                                                &                                                                                                                         &                                                                                                                          & \iddots \arrow[ru, "\nu_{d_{L-2}}"] \arrow[ddd, dashed]         &                                                                                                                                        &                                                                                           \\
                                          &                                                                                                                                                &                                                                                                                         & \mathcal{V}_3\subset\mathbb{P}^{M_3} \arrow[ru, "\nu_{d_4}"] \arrow[dd, "\pi_3", dashed]                                 &                                                                 &                                                                                                                                        &                                                                                           \\
                                          &                                                                                                                                                & \mathcal{V}_2\subset\mathbb{P}^{M_2} \arrow[ru, "\nu_{d_3}"] \arrow[d, "\pi_2", dashed]                                 &                                                                                                                          &                                                                 &                                                                                                                                        &                                                                                           \\
\mathbb{P}^{n_0-1} \arrow[r, "\nu_{d_1}"] & \mathcal{V}^{n_0-1}_{d_1}\subset\mathbb{P}^{N_1} \arrow[ru, "\nu_{d_2}"] \arrow[r, "\overline{\nu}_{d_2}"] \arrow[rd, "\widetilde{\nu}_{d_2}"] & \overline{\mathcal{V}}_2\subset\mathbb{P}^{R_2} \arrow[r, "\overline{\nu}_{d_3}"] \arrow[d, "\overline{\pi}_2", dashed] & \overline{\mathcal{V}}_3\subset\mathbb{P}^{R_3} \arrow[r, "\overline{\nu}_{d_4}"] \arrow[dd, "\overline{\pi}_3", dashed] & \dots \arrow[r, "\overline{\nu}_{d_{L-2}}"] \arrow[ddd, dashed] & \overline{\mathcal{V}}_{L-2}\subset\mathbb{P}^{R_{L-2}} \arrow[r, "\overline{\nu}_{d_{L-1}}"] \arrow[dddd, "\overline{\pi}_{L-2}", dashed] & \overline{\mathcal{V}}_{L-1}\subset\mathbb{P}^{R_{L-1}} \arrow[ddddd, "\overline{\pi}_{L-1}", dashed] \\
                                          &                                                                                                                                                & \widetilde{\mathcal{V}}_2\subset\mathbb{P}^{N_2} \arrow[rd, "\widetilde{\nu}_{d_3}"]                                    &                                                                                                                          &                                                                 &                                                                                                                                        &                                                                                           \\
                                          &                                                                                                                                                &                                                                                                                         & \widetilde{\mathcal{V}}_3\subset\mathbb{P}^{N_3} \arrow[rd, "\widetilde{\nu}_{d_4}"]                                     &                                                                 &                                                                                                                                        &                                                                                           \\
                                          &                                                                                                                                                &                                                                                                                         &                                                                                                                          & \ddots \arrow[rd, "\widetilde{\nu}_{d_{L-2}}"]                  &                                                                                                                                        &                                                                                           \\
                                          &                                                                                                                                                &                                                                                                                         &                                                                                                                          &                                                                 & \widetilde{\mathcal{V}}_{L-2}\subset\mathbb{P}^{N_{L-2}} \arrow[rd, "\widetilde{\nu}_{d_{L-1}}"]                                           &                                                                                           \\
                                          &                                                                                                                                                &                                                                                                                         &                                                                                                                          &                                                                 &                                                                                                                                        & \widetilde{\mathcal{V}}_{L-1}\subset\mathbb{P}^{N_{L-1}}                                         
\end{tikzcd}
$$
\end{footnotesize}
We explain the notation. Here, $\mathcal{V}_{d_1}^{n_0-1}$ is the image of the degree $d_1$ Veronese embedding of $\mathbb{P}^{n_0-1}$, $\mathcal{V}_2 = \nu_{d_2}(\mathcal{V}_{d_1}^{n_0-1})$, $\mathcal{V}_{3} = \nu_{d_3}(\mathcal{V}_2)$, and in general $\mathcal{V}_{k} = \nu_{d_k}(\mathcal{V}_{k-1})$. 

The map $\overline{\pi}_i\circ\pi_i$ is the linear projection from the space of degree $d_i$ homogeneous polynomials on $\mathbb{P}^{M_{i-1}}$ vanishing along $\mathcal{V}_{i-1}$, where we set $\mathcal{V}_1 = \mathcal{V}_{d_1}^{n_0-1}$.

The map $\pi_i$ is the projection from the space of degree $d_i$ homogeneous polynomials on $\mathbb{P}^{M_{i-1}}$ vanishing along $\nu_{d_{i-1}}(\mathbb{P}^{M_{i-2}})$, and $\overline{\mathcal{V}}_i$ is the image of $\mathcal{V}_i$ via $\pi_i$. 

The maps $\overline{\nu}_{d_i},\widetilde{\nu}_{d_i}$ are the degree $d_i$ Veronese embeddings of $\mathbb{P}^{R_{i-1}}$ and $\mathbb{P}^{N_{i-1}}$ respectively, and $\overline{\pi}_i$ is the linear projection from the space of degree $d_i$ homogeneous polynomials on $\mathbb{P}^{R_{i-1}}$ vanishing along $\overline{\nu}_{d_{i-1}}(\overline{\mathcal{V}}_{d_{i-2}})$.

Set
$$
\nu=\nu_{d_{L-1}}\circ\nu_{d_{L-2}}\circ\dots\circ\nu_{d_1},\qquad
\overline{\nu}=\overline{\nu}_{d_{L-1}}\circ\overline{\nu}_{d_{L-2}}\circ\dots\circ\overline{\nu}_{d_1},\qquad
\widetilde{\nu}=\widetilde{\nu}_{d_{L-1}}\circ\widetilde{\nu}_{d_{L-2}}\circ\dots\circ\widetilde{\nu}_{d_1}.
$$
\end{Construction}

Consider the variety
\renewcommand{\arraystretch}{0.5}
\begin{footnotesize}
$$
\mathcal{X}_{\mathbf n,\mathbf d} =
\overline{
  \begin{array}{l}
    \underset{%
      \substack{p_1^1,\dots,p_{n_1}^1 \in \nu_{d_1}(\mathbb{P}^{n_0-1}) \\ \text{general}}%
    }{\scalebox{2.5}{$\bigcup$}}
    \;
    \underset{%
      \substack{p_1^2,\dots,p_{n_2}^2 \in \nu_{d_2}(\langle p_1^1,\dots,p_{n_1}^1\rangle) \\ \text{general}}%
    }{\scalebox{2.0}{$\bigcup$}}
    \cdots
    
    \underset{%
      \substack{p_1^{L-1},\dots,p_{n_{L-1}}^{L-1} \in \nu_{d_{L-1}}(\langle p_1^{L-2},\dots,p_{n_{L-2}}^{L-2}\rangle) \\ \text{general}}%
    }{\scalebox{1.5}{$\bigcup$}}
    \left\langle p_1^{L-1},\dots,p_{n_{L-1}}^{L-1} \right\rangle 
  \end{array}
}
$$
\end{footnotesize}
in $\mathbb{P}^{M_{L-1}}$, and note that 
$$
\mathcal{V}_{\mathbf n,\mathbf d} = \overline{(\overline{\pi}_{L-1}\circ\pi_{L-1})(\mathcal{X}_{\mathbf n,\mathbf d})} \subset\mathbb{P}^{N_{L-1}}.
$$

\begin{Remark}\label{Rdef}
Set
$$
\mathcal{S}_{\mathbf n,\mathbf d} = \left\lbrace (p_1^j,\dots,p_{n_j}^j) \: | \: p_i^j \in \widetilde{\nu}_{d_j}\left(\left\langle p_1^{j-1},\dots,p_{n_{j-1}}^{j-1} \right\rangle\right),\: j = 1,\dots,L-1;\: i = 1,\dots,n_j\right\rbrace\subset \bigtimes_{j=1}^{L-1} (\mathbb{P}^{N_j})^{n_j},
$$
and consider the incidence variety
$$
\begin{tikzcd}
                                    & {\mathcal{I}_{\mathbf n,\mathbf d} = \left\lbrace ((p_1^j,\dots,p_{n_j}^j),p) \: | \: p \in \left\langle p_1^{L-1},\dots,p_{n_{L-1}}^{L-1}\right\rangle \right\rbrace\subset \mathcal{S}_{\mathbf n,\mathbf d}\times \mathbb{P}^{N_{L-1}}} \arrow[ld, "\varphi"'] \arrow[rd, "\psi"] &                      \\
{\mathcal{S}_{\mathbf n,\mathbf d}} &                                                                                                                                                                                                                                                                                      & \mathbb{P}^{N_{L-1}}
\end{tikzcd}
$$
We have
$$
\dim \mathcal S_{\mathbf n,\mathbf d}
=\sum_{i=1}^{L-1}n_i(n_{i-1}-1),
\qquad
\dim \mathcal I_{\mathbf n,\mathbf d}
=\sum_{i=1}^{L-1}n_i(n_{i-1}-1)+(n_{L-1}-1).
$$
Set
$P:=\sum_{i=1}^{L-1}n_i(n_{i-1}-1)+(n_{L-1}-1)$
and
$S:=\sum_{i=1}^{L-2}n_i(n_{i-1}-1) +\binom{n_{L-2}-1+d_{L-1}}{n_{L-2}-1}-1$.
If
$P\le \min\{S,N_{L-1}\}$,
then the refined expected dimension in Definition~\ref{expdimnL1} is $P$. In this
range $\mathcal V_{\mathbf n,\mathbf d}=\psi(\mathcal I_{\mathbf n,\mathbf d})$
has the refined expected dimension if and only if $\psi$ has finite general fiber.
Hence, if $\mathcal V_{\mathbf n,\mathbf d}$ has the refined expected dimension, then
for a general point $p\in\mathcal V_{\mathbf n,\mathbf d}$ there are only
finitely many compatible nested configurations, and in particular only
finitely many last secant spaces
$\langle p_1^{L-1},\dots,p_{n_{L-1}}^{L-1}\rangle$.

Hence the dimension problem can be compared with the classical theory of secant varieties.
\end{Remark}

The previous construction yields a refined expected dimension when $n_L=1$.

\begin{Definition}[Refined expected dimension for $n_L = 1$]\label{expdimnL1}
When $n_L = 1$ the \textit{refined expected dimension} of $\mathcal V_{\mathbf n,\mathbf d}$ is
$$
\expdim_{\mathrm{ref}}(\mathcal V_{\mathbf n,\mathbf d}) = \min\left\lbrace\sum_{i=1}^L n_i(n_{i-1}-1), \sum_{i=1}^{L-2} n_i(n_{i-1}-1) + \binom{n_{L-2}-1+d_{L-1}}{n_{L-2}-1}-1 , N\right\rbrace
$$
where $N = n_{L} \binom{n_0-1+d}{n_0-1}-n_L$. For $n_L=1$ we use this refinement in place of the expected dimension above, and we say that $\mathcal V_{\mathbf n,\mathbf d}$ is \textit{defective} if $\dim(\mathcal V_{\mathbf n,\mathbf d}) < \expdim_{\mathrm{ref}}(\mathcal V_{\mathbf n,\mathbf d})$.
\end{Definition}

\begin{Remark}
The new piece $\sum_{i=1}^{L-2} n_i(n_{i-1}-1) +
\binom{n_{L-2}-1+d_{L-1}}{n_{L-2}-1}-1$ takes into account the cases in
which the last secant variety fills the projective space spanned by the
last Veronese but the neurovariety does not fill the whole ambient projective space. The $-1$ is the projectivization of that last span. Under the degree assumptions of Theorem~\ref{exdNR}, these additional bounds do not cut the parameter count, as shown in the proof of that theorem.
\end{Remark}

Permutation and rescaling are already taken into account by the projective secant construction. Hence Definition~\ref{def:ident} can be reformulated as follows.

\begin{Corollary}\label{cor:ident}
The polynomial neural network is globally identifiable if and only if, for a general element $F\in\mathcal V_{\mathbf n,\mathbf d}$, there is a unique compatible nested configuration
$(p_1^j,\ldots,p_{n_j}^j)_{j=1}^{L-1}\in\mathcal S_{\mathbf n,\mathbf d}$,
up to permutations at each hidden level, such that
$F\in\langle p_1^{L-1},\ldots,p_{n_{L-1}}^{L-1}\rangle$.
In particular, uniqueness of the last set of points alone is sufficient only when the compatible lower-level secant configurations are themselves uniquely determined.
\end{Corollary}
\begin{Remark}
This is the usual notion of identifiability in tensor decomposition and projective geometry; confront also \cite{AGHKT12} for an application to machine learning. Its relation with non defectiveness is studied for instance in \cite{CM23,MM24}. In Section~\ref{sec:nondef} we will use Theorem~\ref{exdNR} to prove identifiability for a class of architectures.
\end{Remark}

\section{The differential of the parameterization}\label{sec:differential}
Fix $\mathbf n=(n_0,\dots,n_L)$, $\mathbf d=(d_1,\dots,d_{L-1})$, with $L\geq2$, and set $D_j=\prod_{t=1}^jd_t$, $D=D_{L-1}$. Write $W_i=(\alpha^i_{r,s})_{1\leq r\leq n_i,\,1\leq s\leq n_{i-1}}$. We work with the full weight matrices; in the affine gauge the fixed entries are simply omitted.

Set $F^{(1)}=W_1x$, that is $F^{(1)}_r=\sum_{s=1}^{n_0}\alpha^1_{r,s}x_{s-1}$, and, for $2\leq j\leq L-1$, set $F^{(j)}=W_jG^{(j-1)}$, where $G^{(j-1)}_s=(F^{(j-1)}_s)^{d_{j-1}}$. We also set $G^{(0)}_s=x_{s-1}$. The outputs are
$F_a=\sum_{s=1}^{n_{L-1}}\alpha^L_{a,s}G^{(L-1)}_s\in\Sym_DK^{n_0},\, a=1,\dots,n_L$.

Assume now that $n_L=1$ and write $F=F_1$. For $j=1,\dots,L-1$ set
$$
C^{(j)}=(\alpha^L_{1,1},\dots,\alpha^L_{1,n_{L-1}})\prod_{t=L-1}^{j+1}\left(\Diag\bigl(d_t(F^{(t)})^{d_t-1}\bigr)W_t\right)\Diag\bigl(d_j(F^{(j)})^{d_j-1}\bigr),
$$
where the empty product is the identity, and write $C^{(j)}=(C^{(j)}_1,\dots,C^{(j)}_{n_j})$. Since $F^{(j)}=W_jG^{(j-1)}$, the chain rule gives $\frac{\partial F}{\partial\alpha^j_{r,s}}=C^{(j)}_rG^{(j-1)}_s$ for $1\leq j\leq L-1$, while $\frac{\partial F}{\partial\alpha^L_{1,s}}=G^{(L-1)}_s$. Hence, after fixing a monomial basis of $\Sym_DK^{n_0}$, the differential has the block form
$$
\left[\{C^{(1)}_rG^{(0)}_s\}_{r,s}\ \Big|\ \{C^{(2)}_rG^{(1)}_s\}_{r,s}\ \Big|\ \cdots\ \Big|\ \{C^{(L-1)}_rG^{(L-2)}_s\}_{r,s}\ \Big|\ \{G^{(L-1)}_s\}_s\right].
$$
We will use the recursive formulas above rather than this block display. The following example illustrates the different contributions.

\begin{Example}\label{tctc}
Take $\mathbf n=(2,2,2,1)$, $\mathbf d=(3,3)$ and the affine chart $a_{1,2}=a_{2,1}=1$. This chart is different from the uniform convention above but is clearly legitimate. Set $L_1=a_{1,1}x+y$, $L_2=x+a_{2,2}y$, $R=b_{1,1}L_1^3+L_2^3$, $S=b_{2,1}L_1^3+L_2^3$, and $F=c_{1,1}R^3+S^3\in\Sym^9K^2$. At $a_{1,1}=a_{2,2}=0$ we have
$F=c_{1,1}(b_{1,1}y^3+x^3)^3+(b_{2,1}y^3+x^3)^3$.
Write $F=[y_0,\dots,y_9]$ in the basis $[x^9,x^8y,\dots,xy^8,y^9]$. Then
$$
y_0=c_{1,1}+1,\quad y_3=3(c_{1,1}b_{1,1}+b_{2,1}),\quad y_6=3(c_{1,1}b_{1,1}^2+b_{2,1}^2),\quad y_9=c_{1,1}b_{1,1}^3+b_{2,1}^3,
$$
and all the other $y_i$ vanish. Set $\Pi=\{y_1=y_2=y_4=y_5=y_7=y_8=0\}\subset\Pbb^9$. On $\{y_0\neq0\}$ set $f_i=y_i/y_0$.

Varying $a_{2,2}$ produces the monomials $x^8y,x^5y^4,x^2y^7$ and
$$
\frac{\partial f}{\partial a_{2,2}}=
\begin{pmatrix}
9\\[2pt]0\\[2pt]0\\[2pt]6\dfrac{y_3}{y_0}\\[6pt]
0\\[2pt]0\\[2pt]3\dfrac{y_6}{y_0}\\[6pt]0\\[2pt]0
\end{pmatrix}.
$$
Similarly, varying $a_{1,1}$ produces $x^7y^2,x^4y^5,xy^8$ and
$$
\frac{\partial f}{\partial a_{1,1}}=
\begin{pmatrix}
0\\[2pt]3\dfrac{y_3}{y_0}\\[6pt]0\\[2pt]0\\[2pt]
6\dfrac{y_6}{y_0}\\[6pt]0\\[2pt]0\\[2pt]
9\dfrac{y_9}{y_0}\\[6pt]0
\end{pmatrix}.
$$
The other three columns are
$$
\frac{\partial f}{\partial b_{1,1}}=
\begin{pmatrix}
0\\0\\ \dfrac{3c_{1,1}}{y_0}\\0\\0\\
\dfrac{6c_{1,1}b_{1,1}}{y_0}\\0\\0\\
\dfrac{3c_{1,1}b_{1,1}^2}{y_0}
\end{pmatrix},\qquad
\frac{\partial f}{\partial b_{2,1}}=
\begin{pmatrix}
0\\0\\ \dfrac{3}{y_0}\\0\\0\\
\dfrac{6b_{2,1}}{y_0}\\0\\0\\
\dfrac{3b_{2,1}^2}{y_0}
\end{pmatrix},
$$
and
$$
\frac{\partial f}{\partial c_{1,1}}=
\begin{pmatrix}
0\\0\\ \dfrac{3(b_{1,1}-b_{2,1})}{y_0^2}\\0\\0\\
\dfrac{3(b_{1,1}^2-b_{2,1}^2)}{y_0^2}\\0\\0\\
\dfrac{b_{1,1}^3-b_{2,1}^3}{y_0^2}
\end{pmatrix}.
$$
The first two columns are independent and have support outside $\Pi$, while the last three are independent for general $(b_{1,1},b_{2,1},c_{1,1})$ and are supported on $\Pi$. Hence $\rank Jf=5$, and $\mathcal V_{(2,2,2,1),(3,3)}$ has the expected dimension.
\end{Example}

\section{Non defectiveness of neurovarieties}\label{sec:nondef}
We now study the dimension of neurovarieties. We begin with a non defective example outside the degree range of the main theorem.
\begin{Example}\label{guid-ex}
Take $L=3$, $(n_0,n_1,n_2,n_3)=(2,3,2,1)$ and
$(d_1,d_2)=(4,3)$. Write
$F=c_{1,1}R^3+S^3$,
where
$R=b_{1,1}L_1^4+b_{1,2}L_2^4+L_3^4, \, S=b_{2,1}L_1^4+b_{2,2}L_2^4+L_3^4$,
$L_1=a_{1,1}x_0+x_1,\, L_2=a_{2,1}x_0+x_1,\, L_3=a_{3,1}x_0+x_1$,
and write $L_3=aL_1+bL_2$.

A direct expansion produces thirteen auxiliary polynomials spanning a $10$-dimensional space, hence they do not give coordinates. We compute the dimension directly in the monomial basis of $\Sym^{12}K^2$. The eight Jacobian columns corresponding to
$a_{1,1},a_{2,1},a_{3,1},b_{1,1},b_{1,2},b_{2,1},b_{2,2},c_{1,1}$ have rank
$8$ for general parameters. The first three describe infinitesimal movements
of $L_1^{12},L_2^{12},L_3^{12}$ on the degree $12$ rational normal curve,
while the last five give the last secant block. Hence
$\dim\mathcal V_{(2,3,2,1),(4,3)}=8$.
Since $d_1=4<2n_1-1=5$, Theorem~\ref{exdNR} does not apply. Hence the degree bound is sufficient but not necessary.
\end{Example}

\subsection*{A linear independence lemma}

We keep the notation of Section~\ref{sec:differential}. Set $D_0=1$ and
$D_j=d_1\cdots d_j$ for $j\ge1$.

\begin{Lemma}\label{lem:independence}
Assume that $n_i\ge2$ for $i=0,\ldots,j$ and $d_i\ge2n_i-1$ for $i=1,\ldots,j$.
Set $A_0:=\langle x_0,\ldots,x_{n_0-1}\rangle$ and, recursively,
$$
A_h:=\sum_{i=1}^{n_h}\bigl(F^{(h)}_i\bigr)^{d_h-1}A_{h-1}\subset\Sym_{D_h}K^{n_0},\qquad h=1,\ldots,j.
$$
Then, for general parameters,
$A_h= \bigoplus_{i=1}^{n_h} \bigl(F^{(h)}_i\bigr)^{d_h-1}A_{h-1}$
for every $h=1,\ldots,j$. In particular,
$\dim A_h=\prod_{t=0}^{h}n_t$.
\end{Lemma}

\begin{proof}
It is enough to exhibit one choice of the parameters for which all the sums
are direct, since the required rank conditions are open.

We use $x_0,x_1$ as distinguished variables. We specialize the parameters
recursively so that, at level $h$,
$F^{(h)}_1=x_0^{D_{h-1}},\, F^{(h)}_2=x_1^{D_{h-1}}$,
and, for $i=3,\ldots,n_h$,
$F^{(h)}_i=x_0^{D_{h-1}}+\lambda_{h,i}x_1^{D_{h-1}}$,
where the $\lambda_{h,i}$ are distinct and non-zero. At the first level we choose the rows of $W_1$ in this way; at every subsequent level this specialization is obtained by using only the first two activated neurons of the preceding level.

We prove the statement by induction on $h$. At every step $A_{h-1}$ has the
form
$A_{h-1}=B_{h-1}\oplus \bigoplus_{r=2}^{n_0-1}x_rC_{h-1}$,
where
$B_{h-1}\subset K[x_0,x_1]_{D_{h-1}},\, C_{h-1}\subset K[x_0,x_1]_{D_{h-1}-1}$,
and
$x_0^{D_{h-1}},x_1^{D_{h-1}}\in B_{h-1}$.
For $h=1$ this follows from
$A_0=\langle x_0,x_1\rangle\oplus
\bigoplus_{r=2}^{n_0-1}Kx_r$.

Write
$X=x_0^{D_{h-1}},\, Y=x_1^{D_{h-1}},\, s=n_h,\, d=d_h$,
and set
$Q_1=X,\, Q_2=Y,\, Q_i=X+\lambda_{h,i}Y\,(i\ge3)$.
Consider a relation
$\sum_{i=1}^s Q_i^{d-1}a_i=0, \, a_i\in A_{h-1}$.
Since the $Q_i$ involve only $x_0,x_1$, the binary part and the parts
containing $x_r$, for $r\ge2$, can be considered separately.

Fix $2\le r\le n_0-1$. If $c_i\in C_{h-1}$ and
$\sum_{i=1}^sQ_i^{d-1}c_i=0$,
write
$c_i=\sum_{\rho=0}^{D_{h-1}-1} c_{i,\rho}x_0^{D_{h-1}-1-\rho}x_1^\rho$.
In the expansion of $Q_i^{d-1}c_i$ the exponent of $x_1$ is of the form
$kD_{h-1}+\rho, \, 0\le k\le d-1,\, 0\le\rho<D_{h-1}$.
This expression determines $k$ and $\rho$ uniquely. Hence, for each fixed
$\rho$, we obtain a relation among
$Q_1^{d-1},\ldots,Q_s^{d-1}$. These forms are linearly independent since
the $Q_i$ are distinct points of $\Pbb^1$ and $d\ge s$. Therefore all the
$c_{i,\rho}$ vanish.

It remains to consider the binary part. Write
$b_i=\alpha_iX+\beta_iY+b_i^{\circ}$,
where $b_i^{\circ}$ contains only monomials
$x_0^{D_{h-1}-\rho}x_1^\rho$ with
$1\le\rho\le D_{h-1}-1$. The same argument shows that all the coefficients
of the $b_i^{\circ}$ vanish independently. We are left with
$\sum_{i=1}^sQ_i^{d-1}(\alpha_iX+\beta_iY)=0$.
The $2s$ forms
$Q_i^{d-1}X,\, Q_i^{d-1}Y, \, i=1,\ldots,s$,
are linearly independent when $d+1\ge2s$. Indeed, after identifying $X,Y$
with coordinates on a two-dimensional vector space, these are the tangent
directions to the degree-$d$ rational normal curve at $s$ distinct points,
and the corresponding interpolation conditions are independent when
$2s\le d+1$. Since $d\ge2s-1$, it follows that
$\alpha_i=\beta_i=0$ for every $i$.

Hence
$\sum_{i=1}^sQ_i^{d-1}A_{h-1}$
is a direct sum. Moreover, since $Q_1=X$, $Q_2=Y$ and $X,Y\in B_{h-1}$,
the new binary part contains
$X^d=x_0^{D_h},\, Y^d=x_1^{D_h}$.
Hence the same description holds for $A_h$, and the induction continues.

Writing the generators obtained recursively from $A_0$ in a fixed monomial
basis, the direct-sum condition is the non-vanishing of a minor. We have
found a specialization where this minor is non-zero, so the result holds for
general parameters. The dimension formula follows inductively from the direct-sum
decomposition and $\dim A_0=n_0$.
\end{proof}

\begin{Remark}\label{rem:sharpness-degree}
The bound $d_i\ge2n_i-1$ is sharp for the argument of Lemma~\ref{lem:independence}, but it is not necessary for non defectiveness. Indeed, when $n_0=2$, the space $A_1$ should have dimension $2n_1$ and is contained in $\Sym_{d_1}K^2$, which has dimension $d_1+1$. Hence the direct sum in Lemma~\ref{lem:independence} can hold only if $d_1\ge2n_1-1$. On the other hand, Example~\ref{guid-ex} is non defective although $d_1=4<2n_1-1=5$. For $n_0>2$ this dimension count no longer forces the bound.

The bound is nevertheless sharp in order of magnitude. Consider $(2,s,r,1)$, with $2\le r<s$, and activation degrees $(d,e)$. Assume $e\ge2r-1$, $de\ge s(r+1)-1$ and $\binom{s-1+e}{s-1}\ge sr$, so that all the hypotheses of Theorem~\ref{exdNR} except possibly $d\ge2s-1$ are satisfied and the ambient bounds do not cut the parameter count $s(r+1)-1$. Set $Y_{d,e}:=\{[G^e]\mid [G]\in\Pbb(\Sym_dK^2)\}$. Since $\dim Y_{d,e}=d$, the neurovariety is contained in $\Sec_r(Y_{d,e})$, and hence $\dim\mathcal V_{(2,s,r,1),(d,e)}\le r(d+1)-1$. Therefore it is defective whenever $r(d+1)<s(r+1)$, that is $d<s+s/r-1$. For $r=2$ this gives defective families for $d<(3s-2)/2$. Thus no sublinear bound in $n_i$ can hold in general. More precisely, the optimal universal asymptotic constant in a sufficient condition of the form $d_i\ge c n_i+O(1)$, if it exists, lies between $3/2$ and $2$.

If one allows a hidden layer of width one, the estimate becomes sharp also in the constant. Indeed, for $(2,s,1,1)$ with activation degrees $(d,e)$ and $e\ge2$ one has $\dim\mathcal V_{(2,s,1,1),(d,e)}=\min\{2s-1,d\}$, so the neurovariety is defective exactly when $d\le2s-2$. This family violates the hypothesis $n_i\ge2$ of Theorem~\ref{exdNR}, and therefore shows the sharpness of the inequality $d_1\ge2n_1-1$ rather than of the theorem as stated. Notice also that the phenomenon is genuinely related to depth: for $(2,s,1)$ the corresponding secant variety of the rational normal curve is non defective. The same dimension count gives a local obstruction at any intermediate level for which $(n_{i-1},n_i,n_{i+1})=(2,s,r)$.
\end{Remark}

\begin{Proposition}\label{prop:differential-rank}
Assume that $n_i\ge2$ for $i=0,\ldots,L-1$ and
$d_i\ge2n_i-1 \,\text{for }i=1,\ldots,L-1$.
At a general parameter point, the only variations in the kernel of the
differential of the projective parameterization $\phi$ are the usual
rescaling directions of the hidden neurons and the projective rescalings of
the outputs.
\end{Proposition}

\begin{proof}
For each $j$ let $A_j$ be the space in Lemma~\ref{lem:independence}. First note that
$F^{(j)}_i\in A_{j-1}$
and that, for every infinitesimal variation of the parameters,
$\dot F^{(j)}_i\in A_{j-1}$.
For $j=1$ this is clear since the $F^{(1)}_i$ are linear forms. For
$j\ge2$ we have
$F^{(j)}_i= \sum_k\alpha^j_{i,k}\bigl(F^{(j-1)}_k\bigr)^{d_{j-1}}$,
while
$$
\dot F^{(j)}_i=
\sum_k\dot\alpha^j_{i,k}\bigl(F^{(j-1)}_k\bigr)^{d_{j-1}}
+d_{j-1}\sum_k\alpha^j_{i,k}
\bigl(F^{(j-1)}_k\bigr)^{d_{j-1}-1}\dot F^{(j-1)}_k.
$$
The assertion follows by induction.

Write the outputs as
$F_a= \sum_{i=1}^{n_{L-1}} \alpha^L_{a,i}\bigl(F^{(L-1)}_i\bigr)^{d_{L-1}}, \, a=1,\ldots,n_L$.
A variation belongs to the kernel of the differential of the map to the
product of projective spaces if, for every $a$, there is a scalar $\mu_a$ such that
$\dot F_a=\mu_aF_a$.
Therefore
$$
0=\sum_{i=1}^{n_{L-1}}
\bigl(F^{(L-1)}_i\bigr)^{d_{L-1}-1}
\left[
(\dot\alpha^L_{a,i}-\mu_a\alpha^L_{a,i})F^{(L-1)}_i
+d_{L-1}\alpha^L_{a,i}\dot F^{(L-1)}_i
\right].
$$
The terms in brackets belong to $A_{L-2}$. By
Lemma~\ref{lem:independence}, the sum is direct. Hence, for every $a,i$,
$(\dot\alpha^L_{a,i}-\mu_a\alpha^L_{a,i})F^{(L-1)}_i +d_{L-1}\alpha^L_{a,i}\dot F^{(L-1)}_i=0$.
For general parameters $\alpha^L_{a,i}\ne0$, and therefore
$\dot F^{(L-1)}_i=\lambda_iF^{(L-1)}_i$
for suitable scalars $\lambda_i$. Moreover,
$\dot\alpha^L_{a,i} =(\mu_a-d_{L-1}\lambda_i)\alpha^L_{a,i}$.
These are exactly the variations obtained by rescaling the $i$-th hidden neuron, compensating in the last matrix, and rescaling the $a$-th output. After subtracting them we may assume
$\dot F^{(L-1)}_i=0,\, \dot W_L=0$.

If $L=2$, the forms $F^{(1)}_i$ are linear and we are already done after
subtracting the rescalings above. Assume $L\ge3$. For every
$u=1,\ldots,n_{L-1}$ we have
$$
0=\dot F^{(L-1)}_u
=\sum_{i=1}^{n_{L-2}}
\bigl(F^{(L-2)}_i\bigr)^{d_{L-2}-1}
\left[
\dot\alpha^{L-1}_{u,i}F^{(L-2)}_i
+d_{L-2}\alpha^{L-1}_{u,i}\dot F^{(L-2)}_i
\right].
$$
Applying Lemma~\ref{lem:independence} at level $L-2$ gives
$\dot F^{(L-2)}_i=\lambda'_iF^{(L-2)}_i$.
These are the usual rescalings at level $L-2$. After subtracting them, $\dot F^{(L-2)}_i=0$. Repeating the argument down to the first level gives
$\dot F^{(1)}_i=0$.
Since the $F^{(1)}_i$ are linear forms in the input variables, the remaining
variations of $W_1$ vanish. Hence the kernel consists exactly of the stated
rescaling directions.
\end{proof}

The following statement is also a consequence of
\cite[Corollary~17 and Proposition~10]{UDBC25}. Indeed, their result gives
finite identifiability, hence maximal affine dimension, under the same degree
bounds. Passing to the product of projective output spaces subtracts the
$n_L$ independent output rescalings and gives the dimension below. We include
the following independent proof since it identifies the kernel of the
differential explicitly.

\begin{thm}\label{exdNR}
Let $\mathbf n=(n_0,\ldots,n_L)$ and
$\mathbf d=(d_1,\ldots,d_{L-1})$, with $L\ge2$. Assume that
$n_i\ge2\,\text{for }i=0,\ldots,L-1$
and
$d_i\ge2n_i-1\,\text{for }i=1,\ldots,L-1$.
Then the neurovariety $\mathcal V_{\mathbf n,\mathbf d}$ is non defective.
More precisely,
$\dim\mathcal V_{\mathbf n,\mathbf d} =\sum_{i=1}^L n_i(n_{i-1}-1)$.
\end{thm}

\begin{proof}
The parameter space has dimension
$\sum_{i=1}^L n_in_{i-1}$.
By Proposition~\ref{prop:differential-rank}, at a general point the kernel of
the differential consists exactly of the usual rescaling directions. There
is one such direction for every neuron at the hidden levels
$1,\ldots,L-1$ and one projective rescaling for every output. These
directions are linearly independent at a general parameter point: looking at
the variation of $W_1$ first forces all the rescaling coefficients at level
$1$ to vanish in any linear relation, and then the same argument applied
successively to $W_2,\ldots,W_{L-1}$ and finally to $W_L$ forces all the
remaining coefficients to vanish. Hence the kernel has dimension
$\sum_{i=1}^L n_i$.
Therefore
$$
\dim\mathcal V_{\mathbf n,\mathbf d} =\sum_{i=1}^L n_in_{i-1}-\sum_{i=1}^L n_i =\sum_{i=1}^L n_i(n_{i-1}-1).
$$

It remains to check that the bounds entering the expected dimension do not cut this number. Set
$D_i=\prod_{h=1}^i d_h,\, B_i=\binom{n_0-1+D_i}{n_0-1}$.
If $m,s\ge2$ and $d\ge2s-1$, then
$\binom{m-1+d}{m-1} \ge1+(m-1)d \ge ms$.
Indeed, the first inequality follows by counting the monomials
$x_1^d$ and $x_1^{d-k}x_r^k$, with $r=2,\ldots,m$ and $k=1,\ldots,d$.
Taking $m=n_0$, $s=n_1$ gives $B_1\ge n_0n_1$. Moreover, for $i\ge2$,
$B_i-B_{i-1}\ge(d_i-1)D_{i-1}\ge(n_i-1)n_{i-1}$.
Hence, by telescoping,
$$
B_{L-1}-1\ge \sum_{i=1}^{L-1}n_i(n_{i-1}-1)+(n_{L-1}-1).
$$
When $n_L=1$ this is the parameter count above. Furthermore, taking
$m=n_{L-2}$ and $s=n_{L-1}$ in the same elementary estimate gives
$$
\binom{n_{L-2}-1+d_{L-1}}{n_{L-2}-1} \ge n_{L-1}n_{L-2},
$$
so the second term in Definition~\ref{expdimnL1} does not cut the parameter
count either. Hence the refined expected dimension of Definition~\ref{expdimnL1}
coincides with the dimension computed above.

For arbitrary $n_L$ the ambient product has dimension
$n_L(B_{L-1}-1)$. Set
$P:=\sum_{i=1}^{L-1}n_i(n_{i-1}-1)$.
Since $B_{L-1}-1\ge P+(n_{L-1}-1)$, we have
$$
n_L(B_{L-1}-1) \ge n_LP+n_L(n_{L-1}-1) \ge P+n_L(n_{L-1}-1) =\sum_{i=1}^L n_i(n_{i-1}-1).
$$
Hence the ambient space does not cut the parameter count, and the result
follows.
\end{proof}

\medskip
For comparison, we record the numerical conditions that detect some possible
failures of the expected differential rank. The room condition is
\stepcounter{thm}
\begin{equation}\label{ineq:room}
n_{j-1}+n_j-1<\binom{n_{j-1}-1+d_j}{n_{j-1}-1}
\qquad\text{for }j=1,\ldots,L-1.
\end{equation}
For $j=1,\ldots,L-2$ set
$$
M_j:=\binom{n_{j-1}-1+d_j}{n_{j-1}-1},\qquad
r_j:=\min\{n_j,n_{j+1}\},\qquad
X_j:=V_{d_j}^{n_{j-1}-1}\subset\Pbb^{M_j-1}.
$$
Let $GS_{X_j}(r_j-1,n_j)\subset\mathbb G(r_j-1,M_j-1)$ be the closure of
the set of $(r_j-1)$-planes contained in the span of $n_j$ general points of
$X_j$. The inequalities
\stepcounter{thm}
\begin{equation}\label{eq:G1}
n_j(n_{j-1}-1)+r_j(n_j-r_j)\le r_j(M_j-r_j)
\end{equation}
and
\stepcounter{thm}
\begin{equation}\label{eq:G2}
\dim GS_{X_j}(r_j-1,n_j)
=n_j(n_{j-1}-1)+r_j(n_j-r_j)
\end{equation}
mean that the natural Grassmann-secant incidence has finite general fiber.
We refer to (\ref{eq:G1})--(\ref{eq:G2}) as the intermediate Grassmann
condition \textup{(G)}. We use this condition below only to describe an
obstruction: the proof of Theorem~\ref{exdNR} does not rely on it.

\begin{Remark}
The room condition alone does not exclude a positive-dimensional family of intermediate secant spaces producing the same higher-level data. The
examples in Remark~\ref{rem:necessity-room} show both the boundary failure of
the room condition and a failure of (\ref{eq:G1}) although the room condition
holds. The latter example was found using the techniques of
\cite{DR2026-MFA-PNN}. The degree assumptions of
Theorem~\ref{exdNR} give the required independence directly at the level of
the differential.
\end{Remark}

\subsection*{Obstructions outside the degree bound}

Example~\ref{guid-ex} is non defective and lies outside the degree range of Theorem~\ref{exdNR}. There the direct sum in Lemma~\ref{lem:independence} already fails at the first level:
$\dim A_1=5<6$. Hence the lemma is sufficient for full differential rank but does not give a characterization. The room and Grassmann-secant counts provide geometric obstructions in this range.

\begin{Remark}\label{rem:necessity-room}
Take $L=3$, $(n_0,n_1,n_2,n_3)=(2,3,2,1)$ and $(d_1,d_2)=(3,3)$.
The hypothesis (\ref{ineq:room}) fails already at the first step:
$n_0+n_1-1 \;=\; 2+3-1 \;=\; 4 \;=\; \binom{n_0-1+d_1}{\,n_0-1\,}$.
Here $M_1=4$ and $r_1=2$, so the source of the Grassmann-secant incidence has dimension
$3\cdot1+2(3-2)=5$, whereas $\dim\mathbb G(1,3)=4$. Hence \textup{(G1)} also fails, with excess one.
We show that, correspondingly, the neurovariety is defective. Write
$$
F \;=\; c_{1,1}\,R^3 + S^3,
\qquad
R \;=\; b_{1,1}L_1^3+b_{1,2}L_2^3+L_3^3,
\qquad
S \;=\; b_{2,1}L_1^3+b_{2,2}L_2^3+L_3^3,
$$
with
$
L_i=a_{i,1}x_0+x_1\ (i=1,2,3)
$.
Hence $L_1,L_2,L_3$ are three points on the twisted cubic
$C=\nu_3(\Pbb^1)\subset \Pbb^3$, and $\nu_3(H)\cong \nu_3(\Pbb^2)$ is the Veronese surface with
$H:=\Span{L_1^3,L_2^3,L_3^3}\subset \Pbb^3$.
The neurovariety is then defective since any plane containing the line $\left\langle R,S\right\rangle$ yields a Veronese surface having $\left\langle R^3,S^3\right\rangle$ as secant line.

The failure of the room condition is accompanied by a degeneration of the adapted evaluation frame. We record an explicit relation among the polynomials occurring in that computation. Consider in $\Sym^9 K^2$ the polynomials
$$
\begin{aligned}
&B_1:=R^2L_1^2x,\quad B_2:=S^2L_1^2x,\quad
B_3:=R^2L_2^2x,\quad B_4:=S^2L_2^2x,\\
&B_7:=R^2L_1^3,\quad B_8:=S^2L_1^3,\quad
B_9:=R^2L_2^3,\quad B_{10}:=S^2L_2^3,\quad
B_{11}:=R^3.
\end{aligned}
$$
Write each $B_k$ in the monomial basis $\{x^9,x^8y,\dots,xy^8,y^9\}$ and denote by
$\mathbf C_k\in K^{10}$ the corresponding coefficient column vector.
Set
$$
A:=a_{1,1}-a_{2,1},\qquad
\begin{array}{lcl}
c_1 & := & 3\,(a_{1,1}-a_{3,1})\,(a_{2,1}-a_{3,1})^{2}\,A,\\
c_3 & := & 3\,(a_{1,1}-a_{3,1})^{2}\,(a_{2,1}-a_{3,1})\,A,\\
c_7 & := & -\,(a_{2,1}-a_{3,1})^{2}\,\bigl(3a_{1,1}-a_{2,1}-2a_{3,1}\bigr)\;-\;b_{1,1}\,A^{3},\\
c_9 & := & -\,(a_{1,1}-a_{3,1})^{2}\,\bigl(a_{1,1}-3a_{2,1}+2a_{3,1}\bigr)\;-\;b_{1,2}\,A^{3},\\
c_{11} & := & A^{3}.
\end{array}
$$
A direct computation shows the non-trivial linear relation
$$
c_1\,B_{1}\;+\;c_3\,B_{3}\;+\;c_7\,B_{7}\;+\;c_9\,B_{9}\;+\;c_{11}\,B_{11}\;=\;0 \quad\text{in }\Sym^9 K^2.
$$
Hence the polynomials
$B_1,B_2,B_3,B_4,B_7,B_8,B_9,B_{10},B_{11}$
are linearly dependent.

This relation is consistent with the positive-dimensional fiber detected above, but by itself it does not determine the rank of the Jacobian. A direct computation of the Jacobian in the monomial basis of $\Sym^9K^2$ gives rank $7$ at a general point, rather than the expected $8$. Hence the neurovariety is defective of defect $1$.

The following example was found using the techniques of \cite{DR2026-MFA-PNN}. Take
$(n_0,n_1,n_2,n_3)=(4,16,10,1),\, (d_1,d_2)=(3,3)$.
Unlike the preceding example, the room condition holds:
$4+16-1=19<\binom{6}{3}=20, \, 16+10-1=25<\binom{18}{15}=816$,
and the last Veronese $V^{15}_3$ is not $10$-defective. Nevertheless, the first intermediate Grassmann condition fails numerically. Indeed, $M_1=20$ and $r_1=10$, hence
$16(4-1)+10(16-10)=108>100=10(20-10)=\dim\mathbb G(9,19)$.
Geometrically, the sixteen points $L_1^3,\ldots,L_{16}^3$ span a general $15$-plane $H\subset\Pbb^{19}$, while the ten second-layer forms span a $9$-plane $\Lambda\subset H$. The family of $16$-secant $15$-planes to $V^3_3$ containing a general such $\Lambda$ therefore has dimension at least
$
108-100=8.
$
A direct Jacobian computation gives dimension $200$ for the affine cone, hence
$\dim\mathcal V_{(4,16,10,1),(3,3)}=199$, while Definition~\ref{expdimnL1} gives refined expected projective dimension $207$. Hence the defect is $8$. This shows that the room condition alone does not guarantee the expected differential rank.
\end{Remark}

\medskip
\noindent\textbf{Uniqueness of the compatible configuration.}
Let $X=V_d^{m-1}$ and assume $d\ge2s-1$. Any set of at most $2s$
distinct points of $\Pbb^{m-1}$ has linearly independent degree-$d$
Veronese images. Indeed, for each point one can separate it from all the
others by a product of at most $2s-1$ hyperplanes and, if necessary,
multiply by a form not vanishing at the chosen points.

It follows that, for every $1\le r\le s$, a general $(r-1)$-plane contained
in the span of $s$ general points of $X$ determines those $s$ points
uniquely. If two sets of $s$ points gave the same general $(r-1)$-plane,
linear independence of the union of the two sets would imply that the
intersection of their spans is the span of their common points. A general
$(r-1)$-plane is not contained in the span of a proper subset of either
set, so the two sets coincide.

Under the assumptions $d_i\ge2n_i-1$, this argument can be applied
successively from the last hidden layer to the first one. Hence a general
compatible last secant space determines a unique compatible nested
configuration in $\mathcal S_{\mathbf n,\mathbf d}$, up to permutations at
the hidden levels.

We can now prove global identifiability for multi-output neurovarieties.
\begin{thm}\label{th:main}
Set an architecture associated to the width vector
${\mathbf n}=(n_0,\ldots,n_L)$ and activation degrees
${\mathbf d}=(d_1,\ldots,d_{L-1})$. Assume that $L\ge2$, $n_L\ge2$,
$n_i\ge2\,\text{for }i=0,\ldots,L-1, \, d_i\ge2n_i-1\,\text{for }i=1,\ldots,L-1$.
Then the associated polynomial neural network is globally identifiable. In
particular, $\mathcal V_{\mathbf n,\mathbf d}$ is non defective.
\end{thm}

\begin{proof}
Non defectiveness follows from Theorem~\ref{exdNR}. Consider the associated
single-output neurovariety
$Y:=\mathcal V_{(n_0,\ldots,n_{L-1},1),\mathbf d}$.
By Theorem~\ref{exdNR},
$\dim Y= \sum_{i=1}^{L-1}n_i(n_{i-1}-1)+(n_{L-1}-1)$.
The proof of Theorem~\ref{exdNR} also shows that the other two terms in the
refined expected dimension of Definition~\ref{expdimnL1} do not cut this number. Hence Remark~\ref{Rdef}
applies, and the incidence map defining $Y$ has finite general fiber. In
particular, through a general point $x_1\in Y$ there are only finitely many
compatible $(n_{L-1}-1)$-planes
$\Lambda=\langle p_1^{L-1},\ldots,p_{n_{L-1}}^{L-1}\rangle$.

Fix one such general $\Lambda$ and choose the outputs
$x_1,\ldots,x_{n_L}$ generally in $\Lambda$. Since $n_{L-1}\ge2$,
$\Lambda$ has positive dimension. If $\Lambda'\ne\Lambda$ is another
compatible plane through $x_1$, then $\Lambda\cap\Lambda'$ is a proper
closed subset of $\Lambda$. Since there are only finitely many such
$\Lambda'$, a general choice of $x_2$ excludes all of them. Hence $\Lambda$
is the unique compatible last secant space containing all the outputs.

By the uniqueness statement above, $\Lambda$ determines the compatible
nested configuration uniquely, up to permutations at the hidden levels.
Corollary~\ref{cor:ident} gives global identifiability.
\end{proof}

\begin{Remark}\label{rem:nsb}
The condition $n_{L-1}\ge2$, which follows from the hypotheses of
Theorem~\ref{th:main}, is used to separate the finitely many compatible last
secant spaces with a second general output. If $n_{L-1}=1$, the common last
secant space is a point, so all projective outputs coincide and this argument
cannot distinguish different compatible configurations. This occurs, for
instance, for the architecture $(3,2,1,2)$ with quadratic activations. The
associated single-output architecture is already outside the non defective
situation considered above: $\Sec_2(V_2^2)\subset\Pbb^5$ is defective.
\end{Remark}

\section{Neurovarieties Magma library}\label{MagmaNV}
A Magma library for computations with neurovarieties is available at
\begin{center}
\url{https://github.com/msslxa/Neurovarieties}
\end{center}
Given $\mathbf n=[n_0,\dots,n_L]$ and $\mathbf d=[d_1,\dots,d_{L-1}]$, the library constructs the symbolic network, the coefficient map, an affine gauge and the Jacobian used to compute the dimension.

\medskip
\noindent\textbf{\texttt{BuildNetworkDataNamed(n, d)}.} This routine returns a polynomial ring $R$, the matrices $W_1,\dots,W_L$ and the output polynomials $v=(F_1,\dots,F_{n_L})$ obtained by formal propagation with $\sigma_i(x)=x^{d_i}$.

\medskip
\noindent\textbf{\texttt{ParametrizationNamed(n, d)}.} It returns the projective coefficient map $\phi:\mathbb A^{\#W}\rightarrow\mathbb P^{m-1}$ obtained by expanding the outputs and collecting the coefficients of the degree $\prod_i d_i$ monomials in the input variables.

\medskip
\noindent\textbf{\texttt{AffineGaugeParametrizationNamed(n, d)}.} The last column of every $W_i$ is set equal to $1$. The routine returns the dehomogenized map $\phi_A:\mathbb A^k\rightarrow\mathbb A^{m-1}$, with readable names for the remaining weights.

\medskip
\noindent\textbf{\texttt{GaugeJacobianRankRandom(n, d)}.} It evaluates the Jacobian of $\phi_A$ at a random integral point avoiding the poles and returns its rank and the sampling point. Repeating the computation gives the generic rank and hence the dimension of the image.

\medskip
\noindent\textbf{\texttt{NeuroVarietyStats(n, d : tries := 10)}.} It returns $\langle\mathrm{expdim},\mathrm{dimActual},\mathrm{fiberDim},\mathrm{defective},P\rangle$, where $\mathrm{dimActual}$ is the maximal Jacobian rank found in the prescribed number of trials, $\mathrm{fiberDim}=k-\mathrm{dimActual}$ and $\mathrm{defective}$ records whether $\mathrm{dimActual}<\mathrm{expdim}$.

\medskip
\noindent\textbf{\texttt{ExhaustiveSmallDefectTest}.} These scripts run the previous routines on bounded families of architectures, apply elementary binomial filters and report $(\mathrm{expdim},\mathrm{dimActual},\mathrm{fiberDim},\mathrm{defective})$.

\medskip
\noindent\textbf{Composite Veronese and projection routines.} They construct $\nu_{d_1}\circ\cdots\circ\nu_{d_j}$, compute the linear relations vanishing on its image and the induced linear projection. They are useful for comparing shallow network parameterizations with classical Veronese geometry.

\bibliographystyle{amsalpha}
\bibliography{Biblio}

\end{document}